\documentclass{elsart}

\usepackage{xypic}
\usepackage{amsmath}
\usepackage{amssymb}

\newcommand{\la}{\longrightarrow}

\newcommand{\N}{\mathbb{N}}
\newcommand{\Z}{\mathbb{Z}}

\DeclareMathOperator{\Hom}{Hom} 
 
 \DeclareMathOperator{\ima}{im}

\DeclareMathOperator{\gldim}{gldim}

\DeclareMathOperator{\holim}{\underleftarrow{\textrm{holim}}}

\newcommand{\A}{\mathcal{A}}
\newcommand{\K}{\mathcal{K}}

\newcommand{\CS}{\mathcal{S}}
\newcommand{\T}{\mathcal{T}}
\newcommand{\CP}{\mathcal{P}}

\newcommand{\Ab}{\mathcal{A}b}
\newcommand{\Inj}{\hbox{\rm Inj-}\!}

\newcommand{\opp}{^\textit{o}}
\newcommand{\Mod}[1]{\hbox{\rm Mod}({#1})}
\newcommand{\md}[1]{\hbox{\rm mod}({#1})}

\newcommand{\Add}[1]{\mathrm{Add}({#1})}
\newcommand{\Prod}[1]{\mathrm{Prod}({#1})}
\newcommand{\Prd}{\mathrm{Prod}}

\newcommand{\Htp}[1]{\mathbf{K}({#1})}
\newcommand{\Cpx}[1]{\mathbf{C}({#1})}

\theoremstyle{plain}
\newtheorem{expl}[thm]{Example}

\begin{document}

\begin{frontmatter}

\title{The dual of Brown representability for homotopy categories
of complexes}
\author[cluj]{George Ciprian Modoi\thanksref{cncs}}
\thanks[cncs]{Research supported by CNCS-UEFISCDI grant PN-II-RU-TE-2011-3-0065} \ead{cmodoi@math.ubbcluj.ro}

\address[cluj]{Faculty of Mathematics and Computer Science, Department of
  Mathematics,  "Babe\c s-Bolyai" University,  1, M. Kog\u alniceanu, 400084 Cluj-Napoca, Romania,}

\begin{abstract}
We call product generator of an additive category a fixed object
satisfying the property that every other object is a direct factor
of a product of copies of it. In this paper we start with an additive category with
products and images, e.g. a module category, and we are concerned
with the homotopy category of complexes with entries in that additive category. 
We prove that the Brown
representability theorem is valid for the dual of the homotopy
category if and only if the initial additive category has a
product generator.
\end{abstract}

\begin{keyword}
Brown representability, Adjoint functors, Triangulated category,
Homotopy category of complexes \MSC 18E30, 16D90, 55U35
\end{keyword}

\date{\today}

\end{frontmatter}

\section*{Introduction}

Brown representability is a key tool in the theory of triangulated
categories. Recall that if $\K$ is a  triangulated category with
products then $\K\opp$ is said to {\em satisfy Brown
representability} if every homological product preserving functor
$F:\K\to\Ab$ is representable. Dually $\K$ satisfies Brown
representability if every cohomological (contravariant) functor
which sends coproducts into products $F:\K\to\Ab$ is
representable. Sometime we call these properties Brown
representability for covariant, respectively contravariant
functors.

The notion of well--generated triangulated category was introduced by Neeman
in his influential book \cite{N}, where it is also shown that Brown
representability holds for triangulated categories of this type.
Prototypes for (algebraic) well--generated triangulated categories
are derived categories and their localizations (see \cite{P}). But
until recently, only a very little was known about Brown representability
for homotopy categories of complexes. Only three papers \cite{CN}, \cite{S} and
\cite{MS} gave some information in this sense, the first for the homotopy category 
over an abelian category without a generator and the other two for homotopy category 
over module categories.
The present work comes to complete the picture started in
\cite{MS}. To be precise, let $R$ be a ring. We denote by
$\Htp{\Mod R}$ the homotopy category of complexes of $R$-modules.
In \cite{MS} it is shown that $\Htp{\Mod R}$ satisfies Brown
representability if and only if $R$ is pure--semisimple. But for
the dual $\Htp{\Mod R}\opp$ only one direction was shown: If
$\Htp{\Mod R}\opp$ satisfies Brown representability then $\Mod R$
must have a product generator. Note that a product generator of an
additive category $\A$ is defined to be an object $G$ with the
property that every object of $\A$ is a direct factor of a product
of copies of $G$. The module category over a pure semi--simple
ring $R$ satisfies the dual property, namely $\Mod R=\Add G$ for
some $G\in\Mod R$, where through $\Add G$ we understand the class
of all direct summands of direct sums of copies of $G$. The main
result in this paper proves the equivalence between the conditions
$\Htp{\Mod R}\opp$ satisfies Brown representability and $\Mod R$
has a product generator. Moreover our approach may be easily
dualized in order to give (a generalization of) results in
\cite{MS} about Brown representability for contravariant functors
defined on homotopy category of complexes.

The problem of Brown representability for covariant functors is
difficult and not completely solved even in the case of
well--generated categories. For that reason the
method used to prove this kind of result deserves perhaps a few
words. In \cite{MP} we proved a generalization of Neeman's variant
of Brown representability for contravariant functors defined on
well--generated triangulated categories. With this aim, we
developed a technique, based on the fact that every object of a
well--generated category is the homotopy colimit of a tower of
objects which is constructed iteratively starting with a set. The
whole construction is analogous to the case of an object of an
abelian category which is filtered by a set (see \cite[Definition
3.1.1]{GT}), but as usual, short exact sequences are replaced by
triangles. Naturally appeared the question if the construction may
be dualized in order to give some information about Brown
representability for covariant functors. Strictly in the setting
of \cite{MP} the answer is probably no, but we adapted here
this method and we observed that if $\Mod R$ has a product
generator, then there is a set of complexes, such that every complex in
$\Htp{\Mod R}$ is cofiltered by that
set; roughly speaking, this means every complex is isomorphic to the
homotopy limit of an inverse tower constructed iteratively starting with that set.

The paper is organized as follows: Section \ref{hell} contains a
new proof of an old (but seemingly not largely known) representability
theorem due to Heller, for functors $F:\K\to\Ab$, where $\K$ is a
triangulated category with products. Some applications to much
recent results are also indicated. Using this, we prove in Section
\ref{cofilt} a new representability theorem, supposing in addition
that every object of $\K$ is cofiltered by a set. Next Section
\ref{brdual} contains the main result of this paper: If we consider an
additive category $\A$ with split
idempotents and products, possessing images or kernels, then
$\Htp\A\opp$ satisfies Brown representability
exactly if $\A$ has a product generator. In particular
we apply this for $\A$ being the module category over a ring $R$,
thus $\Htp{\Mod{R}}\opp$ satisfies Brown
representability if and only if $\Mod{R}$ has a product generator.

\section{A new proof for Heller's representability theorem}\label{hell}
Consider a preadditive category $\K$.  We write $\K(K,K')$ for the
abelian group of morphisms between $K$ and $K'$ in $\K$. By a {\em
(right) $\K$-module\/} we understand a functor $X:\K\opp\to\Ab$.
In this paper modules will always be at right, so for dealing with a {\em left $\K$-module}
we have to consider a right $\K\opp$-module, that is a functor
$X:\K\to\Ab$. A $\K$-module is called {\em finitely presentable\/}
if there is an exact sequence of functors
\[\K(-,K_1)\to\K(-,K_0)\to X\to0\] for some $K_0,K_1\in\K$. We denote
$\Hom_\K(X,Y)$ the class of all natural transformations between
two $\K$-modules. Generally there is no reason for this class to
be a set. However, using Yoneda lemma, we know that
$\Hom_\K(X,Y)$ is actually a set, provided that $X$ is finitely
presentable. We consider the category $\md\K$ of all finitely
presentable $\K$-modules, having $\Hom_\K(X,Y)$ as morphisms
spaces, that is $\md\K(X,Y)=\Hom_\K(X,Y)$ for all $X,Y\in\md\K$.

The Yoneda functor
\[H=H_\K:\K\to\md{\K\opp}\opp\hbox{ given by }H_\K(K)=\K(K,-)\] is  an
embedding of $\K$ into $\md{\K\opp}\opp$, according to Yoneda
lemma. Moreover $\md{\K\opp}\opp$ has kernels. If, in addition,
$\K$ has products then $\md{\K\opp}\opp$
is complete and the Yoneda embedding preserves products. It is
also well--known (and easy to prove) that, if $F:\K\to\A$ is a
functor into an additive category with kernels, then there is a
unique, up to a natural isomorphism, kernel preserving functor
$F^*:\md{\K\opp}\opp\to\A$, such that $F=F^*H_\K$ (see \cite[Lemma
A.1]{KL}). Moreover, $F$ preserves products if and only if $F^*$
preserves limits.

Let $F:\K\to\Ab$ be a functor. The {\em category of elements of
$F$}, denoted by $\K/F$, has as objects pairs of the form $(X,x)$
where $X\in\K$ and $x\in F(X)$, and a map between $(X,x)$ and
$(Y,y)$ in $\K/F$ is a map $f:X\to Y$ in $\K$ such that
$F(f)(x)=y$. Recall that the solution set condition for functors
with values in the category of abelian groups $F:\K\to\Ab$ may be
stated as follows: There is a set $\CS$ of objects in $\K$,
such that for any $K\in\K$ and any $y\in\ F(K)$ there are
$S\in\CS$, $x\in F(S)$ and $f:S\to K$ satisfying $F(f)(x)=y$ (see
\cite[Chapter V, \S 6, Theorem 3]{MW}). We may reformulate this by
saying that the category
\[\CS/F=\{(S,x)\mid S\in\CS, x\in F(S)\}\]
is {\em weakly initial} in $\K/F$, that is for every
$(K,y)\in\K/F$ there exists a map $(S,x)\to(K,y)$ for some
$(S,x)\in\CS/F$. Via Yoneda lemma, every object $(S,x)\in\CS/F$
corresponds to a natural transformation $\K(S,-)\to F$. In these
terms, the existence of a solution set is further equivalent to
the fact that there are objects $S_i\in\K$ indexed over a set $I$
and a functorial epimorphism
\[\bigoplus_{i\in I}\K(S_i,-)\to F\to 0.\] We say that $F$ has a
solution object provided that there is an object $S\in\K$ and a
functorial epimorphism \[\K(S,-)\to F\to 0,\] or equivalently, the
category $\K/F$ has a weakly initial object. Note that if there
are arbitrary products in $\K$, and the functor $F$ preserves
them, then the existence of a solution set is clearly equivalent
to that of a solution object. Obviously if $F\cong\K(S,-)$ is
representable, then $F$ has a solution object.

In the rest of this Section the category $\K$ will be triangulated with
split idempotents. For definition and basic properties of
triangulated categories the standard reference is \cite{N}. Note
that $\K$ has split idempotents, provided that $\K$ has
countable coproducts or products, according to \cite[Proposition
1.6.8]{N} or its dual. Recall that $\K$ is supposed to be
additive. A functor $\K\to\A$ into an abelian category $\A$ is
called {\em homological\/} if it sends triangles into exact
sequences. A contravariant functor $\K\to\A$ which is homological
regarded as a functor $\K\opp\to\A$ is called {\em
cohomological\/} (see \cite[Definition 1.1.7 and Remark
1.1.9]{N}). An example of a homological functor is the Yoneda
embedding $H_\K:\K\to\md{\K\opp}\opp$. We know that in this case
$\md{\K\opp}\opp$ is equivalent to $\md\K$ (see \cite[Remark 5.1.19 and what follows]{N}).
Moreover it
is an abelian category, and for every functor $F:\K\to\A$  into an
abelian category, the unique left exact functor
$F^*:\md{\K\opp}\opp\to\A$ extending $F$ is exact if and only if
$F$ is homological, by the dual of \cite[Lemma 2.1]{KS}. Note that
this is the reason for which $\md{\K\opp}\opp$ (or often the
equivalent category $\md\K$) is called the {\em abelianization} of
the triangulated category $\K$. By \cite[Corollary 5.1.23]{N},
$\md{\K\opp}\opp$ is a Frobenius abelian category, with enough
injectives and enough projectives, which are, up to isomorphism,
exactly objects of the form $\K(K,-)$ for some $K\in\K$.

Observe that in the particular case when the codomain of the
homological functor $F$ is the category $\Ab$ of all abelian
groups, then it may be easily seen that
$F^*(X)\cong\Hom_{\K\opp}(X,F)$, naturally for all
$X\in\md{\K\opp}\opp$. Thus we obtain:

\begin{lem}\label{fstar} If $\K$ is a triangulated category
with split idempotents, then a homological functor
$F:\K\to\Ab$ is representable if and only if its extension
$F^*:\md{\K\opp}\opp\to\Ab$ is representable.
\end{lem}

\begin{pf}
As before $F^*(X)\cong\Hom_{\K\opp}(X,F)$, for all $X\in\md{\K\opp}\opp$.
If $F$ is representable, then $F\in\md{\K\opp}\opp$, so
$F^*$ is represented by $F$. Conversely if $F^*$ is representable
by an object in $\md{\K\opp}\opp$ then this object must be
isomorphic to $F$, therefore $F\in\md{\K\opp}\opp$. Because $F^*$
is exact, $F$ must be projective, hence representable (see \cite[Lemma 5.1.11]{N}).
\end{pf}

\begin{lem}\label{solution}
If $\K$ is a triangulated category with split idempotents,
then a cohomological functor $F:\K\to\Ab$ has a solution object if
and only if $F^*:\md{\K\opp}\opp\to\Ab$ has a solution object.
\end{lem}

\begin{pf} Suppose $F$ has a solution object, i.e. there is a
functorial epimorphism  $H(K)=\K(K,-)\to F\to 0$, with $K\in\K$.
In order to show that $F^*$ has a solution object, it is enough to prove
that the induced natural transformation
\[\Hom_{\K\opp}(-,H(K))\to\Hom_{\K\opp}(-,F)\cong F^*\] is an epimorphism. That is, we want
to show that the map
\[\Hom_{\K\opp}(X,H(K))\to\Hom_{\K\opp}(X,F)\] is surjective, for
all $X\in\md{\K\opp}\opp$.  According to \cite[5.1.23]{N} every finitely presentable
$\K\opp$--module $X$ admits
an embedding $0\to X\to H(U)$, that is an epimorphism from the
projective object $H(U)$ to $X$ in the opposite category $\md{\K\opp}\opp$. Since
$H(K)\in\md{\K\opp}\opp$ is projective--injective and $F^*$ is exact, we
obtain a diagram with exact rows:
\[\diagram\Hom_{\K\opp}(H(U),H(K))\rto\dto&\Hom_{\K\opp}(X,H(K))\rto\dto&0\\
            \Hom_{\K\opp}(H(U),F)\rto&\Hom_{\K\opp}(X,F)\rto&0
\enddiagram.\]
By Yoneda lemma we know that the first vertical map is isomorphic
to $\K(K,U)\to F(U)$, hence it is surjective, thus the diagram
above proves the direct implication.

Conversely if there is $X\in\md{\K\opp}\opp$ and a natural
epimorphism
\[\Hom_{\K\opp}(-,X)\to\Hom_{\K\opp}(-,F)\to0,\] then let $H(K)\to X\to 0$
be an epimorphism in $\md{\K\opp}$ (that is a monomorphism in the
opposite direction in $\md{\K\opp}\opp$), with $K\in\K$. Consider
the composed map
\[\Hom_{\K\opp}(-,H(K))\to\Hom_{\K\opp}(-,X)\to\Hom_{\K\opp}(-,F).\]
Evaluating it at $H(U)$ for an arbitrary $U\in\K$, we obtain a
surjective natural map $\K(K,U)\to F(U)$, hence $F$ has a
solution object.
\end{pf}

\begin{thm}\label{freydstyle}\cite[Theorem 1.4]{H}
If $\K$ is a triangulated category with products, then a
homological products preserving functor $F:\K\to\Ab$ is
representable if and only if it has a solution object.
\end{thm}

\begin{pf}
Under the hypotheses imposed on $\K$ and $F$, the abelian category
$\md{\K\opp}\opp$ is complete and the induced functor
$F^*:\md{\K\opp}\opp\to\Ab$ preserves limits. Therefore it is
representable if and only if it has a solution object, by Freyd's
Adjoint Functor Theorem. Thus the conclusion follows by combining
Lemmas \ref{fstar} and \ref{solution}.
\end{pf}

\begin{rem}\label{isold}
{\rm Theorem \ref{freydstyle} says more than the Neeman's Freyd style
representability theorem \cite[Theorem 1.3]{NR}. Indeed the cited
result states that if every cohomological functor which sends
coproducts into products has a solution objects, then every such a
functor is representable, whereas our result involves a fixed
functor. However the result is known: It already appeared in
Heller's paper \cite{H}. We have just proved the dual version
because our argument is different from Heller's one, and it shows
us explicitly how the result follows from  Freyd's celebrated
Adjoint Functor Theorem.}
\end{rem}

In the same sense in which Theorem \ref{freydstyle} above is
an improvement of \cite[Theorem 1.3]{NR}, we may improve \cite[Theorem 3.7]{MP},
which is the main result there and which uses Neeman's result \cite[Theorem 1.3]{NR}
(for the unexplained terms see \cite{MP}):

\begin{cor}
Let $\K$ be a triangulated category with coproducts which is $\aleph_1$-perfectly generated
by a projective class $\CP$. If $F:\K\to\Ab$ is a cohomological functor which sends coproducts
to products, such that $\CP^{*n}/F$ has a weak terminal object for all $n\in\N$, then $F$ is
representable.
\end{cor}

 In a particular case, namely in the presence of products, we may derive from the above results the dual
 of \cite[Proposition 1.4]{NA}.
In order to state this, recall that if $\K$ is a full subcategory
of $\T$ then a {\em$\K$--preenvelope} of $T\in\T$ is a morphism
$T\to X_T$ with $X_T\in\K$ such that the induced map
$\T(X_T,X)\to\T(T,X)$ is surjective for all $X\in\K$. Dually we
define the concept of {\em precover}. The subcategory $\K$ is
called preenveloping is every object in $\T$ admits a
$\K$-preenvelope.

\begin{cor}\label{envelops}
Let $\T$ be a triangulated category with products, and let $\K$ be a colocalizing subcategory.
The following are equivalent:
\begin{itemize}
 \item[{\rm (i)}] The inclusion $\K\to\T$ has a left adjoint.
 \item[{\rm (ii)}] Every object in $\T$ admits a $\K$--preenvelope.
\end{itemize}
\end{cor}

\begin{pf} Since the implication (i)$\Rightarrow$(ii) follows from the general theory
of adjoint functors, we only need to show the converse. But this follows immediately from
Theorem \ref{freydstyle} since, if $I:\K\to\T$ is the inclusion functor, then for every $T\in\T$ the functor
$\T(T,I(-)):\K\to\Ab$ is homological,
preserves products and has a solution object, given by the functorial epimorphism $\K(X_T,-)\to\T(T,I(-))$,
where $T\to X_T$ is a $\K$--preenvelope of $X$.
\end{pf}

\section{Cofiltered objects in triangulated categories}\label{cofilt}

As before we denote by $\K$ a triangulated category with products.
Let $\CS\subseteq\K$ be a set of objects. We denote $\Prod\CS$ the full
subcategory of $\K$ consisting of all direct factors of products of objects in $\CS$.
We define inductively $\Prd_0(\CS)=\Prod\CS$ and $\Prd_n(\CS)$ is
the full subcategory of $\K$ which consists of all objects $Y$
lying in a triangle $X\to Y\to Z\to X[1]$ with $X\in\Prd_0(\CS)$
and $Y\in\Prd_n(\CS)$. We suppose that $\CS$ is closed under
suspensions and desuspensions, so the same is true for
$\Prd_n(\CS)$, by \cite[Remark 07]{NR}. Moreover the same
\cite[Remark 07]{NR} tells us that if $X\to Y\to Z\to X[1]$ is a
triangle with $X\in\Prd_n(\CS)$ and $Y\in\Prd_m(\CS)$ then
$Z\in\Prd_{n+m}(\CS)$. We say that an object $X\in\K$ is
{\em$\CS$-cofiltered} if it may be written as a homotopy limit
$X\cong\holim X_n$ of an inverse tower, with $X_0\in\Prd_0(\CS)$, and
$X_{n+1}$ lying in a triangle $P_{n}\to X_{n+1}\to X_n\to P_n[1],$
for some $P_n\in\Prd_0(\CS)$. Inductively we have
$X_n\in\Prd_n(\CS)$, for all $n\in\N$.

\begin{lem}\label{weaki}
Let $\K$ be a triangulated category and let $\CS\subseteq\K$ be a
set closed under suspensions and desuspensions.
Suppose that every $X\in\K$ is $\CS$-cofiltered. Then every
homological product preserving functor $F:\K\to\Ab$ has a solution
object.
\end{lem}

\begin{pf} We shall prove a statement equivalent to the conclusion, namely
that the category of elements $\T/F$ has a weakly initial object.
In order to do this, we shall apply the dual of the argument used in
the proof of \cite[Proposition 3.6]{MP}. Since the hypotheses are slightly modified,
we sketch here this argument (in the dual form appropriate
to the present approach).

By \cite[Lemma 2.3]{NR}, we know that the category $\Prd_n(\CS)/F$ has a weakly initial
object denoted $(T_n,t_n)$, for all $n\in\N$.
 Let $I$ be the non--empty set of all inverse towers
of the form
\[T_0\stackrel{w_0}\longleftarrow T_1\stackrel{w_1}\longleftarrow T_2\longleftarrow \cdots\]
with $F(w_n)(t_{n+1})=t_n$, for all $n\in\N$, and denote by $T(i)$ the homotopy limit of the
tower $i\in I$. By
\cite[Lemma 5.8(2)]{B}, there is an exact sequence
\[0\to{\underleftarrow\lim}^{(1)}F(T_n[-1])\to F(\holim T_n)\to{\underleftarrow\lim} F(T_n)\to 0.\]
Clearly $(t_n)_{n\in\N}\in{\underleftarrow\lim} F(T_n)$, thus
there exists $t(i)\in F(T(i))=F(\holim T_n)$ which maps in $(t_n)_{n\in\N}$
via the surjective morphism above. Putting $T=\prod_{i\in I}T(i)$ and $t=(t(i))_{i\in I}$ we
claim that $(T,t)$ is a weakly initial object in $\K/F$. In order to prove the claim,
consider an object $X\in\K$. By hypothesis, there is an inverse tower
\[X_0\stackrel{u_0}\longleftarrow X_1\stackrel{u_1}\longleftarrow X_2\longleftarrow \cdots\]
whose homotopy limit is $X$ such that $X_0\in\Prd_0(\CS)$, and
every $X_{n+1}$ lies in a triangle $P_{n}\to
X_{n+1}\stackrel{u_n}\longrightarrow X_{n}\to P_n[1],$ for some
$P_n\in\Prd_0(\CS)$. We use again \cite[Lemma 5.8(2)]{B} for
constructing the commutative diagram with exact rows:
\[\diagram
0\rto&{\underleftarrow\lim}^{(1)}\K(T,X_n[-1])\rto\dto&\K(T,\holim X_n)\rto\dto&\underleftarrow\lim\K(T,X_n)\rto\dto&0\\
0\rto&{\underleftarrow\lim}^{(1)}F(X_n[-1])\rto&F(\holim X_n)\rto&\underleftarrow\lim F(X_n)\rto&0
\enddiagram\]
whose columns are induced by the natural transformations which correspond to $t\in F(T)$
under Yoneda Lemma. If we show that the two extreme vertical arrows are surjective,
the same is true for the middle arrow too, and we are done. But for the first
vertical map this follows by the commutative diagram: \[\diagram
\prod_{n\in\N}\K(T,X_n[-1])\rto\dto&{\underleftarrow\lim}^{(1)}\K(T,X_n[-1])\dto\\
\prod_{n\in\N}F(X_n[-1])\rto&{\underleftarrow\lim}^{(1)}F(X_n[-1])
\enddiagram\] whose arrows connected with the south-west corner are both surjective.

In order to prove that the third vertical map above is surjective,
we consider an element $x\in\underleftarrow\lim F(X_n)$, that is
$x=(x_n)_{n\in\N}\in\prod F(X_n)$ such that $F(u_n)(x_{n+1})=x_n$,
for all $n\in\N$. Next we construct a commutative diagram
\[\diagram
T_0\dto_{f_0}&T_1\lto_{w_0}\dto_{f_1}&T_2\lto_{w_1}\dto^{f_2}&\cdots\lto\\
X_0&X_1\lto^{u_0}&X_2\lto^{u_1}&\cdots\lto
\enddiagram\] whose upper line is a tower in $I$, and satisfying $F(f_n)(t_n)=x_n$ for
all $n\in\N$. This construction is performed inductively as
follows: $f_0$ comes from the fact that $(T_0,t_0)$ is weakly
initial in $\Prd_0(\CS)/F$. Suppose the first $n$ steps are done.
We construct the following commutative diagram whose rows are
triangles and the middle square is homotopy pull-back (see
\cite[Definition 1.4.1]{N}):
\[\diagram
P_n\rto\ddouble&Y_{n+1}\rto\dto&T_n\rto\dto^{f_n}&P_n[1]\ddouble\\
P_n\rto&X_{n+1}\rto^{u_n}&X_n\rto&P_n[1]
\enddiagram\]
The upper triangle shows that $Y_{n+1}\in\Prd_{n+1}(\CS)$ where $(T_{n+1},t_{n+1})$ is
weakly initial, hence we find a map
$(T_{n+1},t_{n+1})\to(Y_{n+1},y_{n+1})$ in $\Prd_{n+1}(\CS)/F$. Now $Y_{n+1}$ is obtained
via a triangle \[Y_{n+1}\to T_n\times X_{n+1}\stackrel{(f_n,-u_n)}\longrightarrow X_n\to Y_{n+1}[1].\]
Applying the homological functor $F$ we get an exact sequence:
\[F(Y_{n+1})\to F(T_n)\times F(X_{n+1})\stackrel{(F(f_n),-F(u_n))}\longrightarrow F(X_n).\]
Since $F(f_n)(t_n)-F(u_n)(x_{n+1})=x_n-x_n=0$ we get an element
$y_{n+1}\in Y_{n+1}$, which maps in $(t_n,x_{n+1})$, via the first
morphism of the exact sequence above.  The morphism $f_{n+1}$ is
the composition  $T_{n+1}\to Y_{n+1}\to X_{n+1}$. The upper row
above is, as noticed, an inverse tower in $I$, and let denote it
by $i$. Finally the element $t\in T$ maps to
$(x_n)_{n\in\N}\in\underleftarrow{\lim}F(X_n)$, via the map
$F(T)\to
F(T(i))\to\underleftarrow{\lim}F(T_n)\to\underleftarrow{\lim}F(X_n)$,
proving that the map
$\underleftarrow{\lim}\K(T,X_n)\to\underleftarrow{\lim}F(X_n)$ is
surjective.
\end{pf}

Combining Theorem \ref{freydstyle} and Lemma \ref{weaki} we obtain:

\begin{thm}\label{cof}
Let $\K$ be a triangulated category. Suppose there is a set $\CS\subseteq\K$ closed
under suspensions and desuspensions, such that
every $X\in\K$ is $\CS$-cofiltered. Then every
homological product preserving functor $F:\K\to\Ab$ is
representable, therefore $\K\opp$ satisfies Brown
representability.
\end{thm}

\section{Brown representability for the dual of a homotopy
category}\label{brdual}

Throughout this section $\A$, will denote an additive category,
that is preadditive, with zero object and
finite biproducts; we suppose also that $\A$ has split idempotents.
We consider categories $\Cpx\A$ and $\Htp\A$ called the {\em category of
complexes} respectively the {\em homotopy category of complexes} over $\A$,
both of them having
as objects complexes of objects in $\A$, that is a chain of objects and morphisms (called {\em differentials})
in $\A$ of the form
\[X=\cdots\to X^{n-1}\stackrel{d_X^{n-1}}\la
X^{n}\stackrel{d_X^{n}}\la X^{n+1}\to\cdots,\]such that $d_X^{n}d_X^{n-1}=0$ for all $n\in\Z$.
The morphisms in the category $\Cpx\A$ are families $(f^n)_{n\in\Z}$ of morphisms in $\A$ commuting with
differentials, and
$\Htp\A(X,Y)=\Cpx\A(X,Y)/\sim$ where $\sim$ is an equivalence relation called {\em homotopy},
defined as follows: two maps of complexes $(f^n)_{n\in\Z},(g^n)_{n\in\Z}:X\to Y$ are homotopically
equivalent if there is $s^n:X^n\to Y^{n-1}$, for all
$n\in\Z$ such that $f^n-g^n=d_Y^{n-1}s^n+s^{n+1}d_X^n$.  Note that $\Cpx\A$ is an exact category (in the sense of \cite[Section 4]{Kel}) with
respect to all
short exact sequences which split in each degree (see \cite[Example 4.3]{Kel}), and
$\Htp\A$ may be constructed as the stable category of this exact category by
\cite[Example 6.1]{Kel}. Hence $\Htp\A$ is a triangulated category. Note that the structure 
of triangulated category comes with a translation functor denoted by $[1]$, where $X[1]^n=X^{n+1}$ and 
$d_{X[1]}^n=-d_X^{n+1}$. 
It is well--known that $\Htp\A$ has (co)products provided that $\A$ does the same.
Considering every object in $\A$ as a complex concentrated in degree $0$, the category $\A$ may be embedded 
in $\Htp\A$. 

Fix the additive category $\A$ as before. For an object $G\in\A$ we denote by $\Prod G$ respectively
$\Add G$ the full subcategory consisting of direct factors (or equivalently, direct summands) of a
product (respectively coproduct) of copies of $G$ (assuming that the requested products or coproducts
exist). We say that $\A$ has a {\em product generator} if there is an object $G\in\A$ such that
$\A=\Prod G$. For the dual situation when $\A=\Add G$ we use the more standard terminology
{\em $\A$ is pure--semisimple} (see \cite[Definition 2.1 and Proposition 2.2]{S}).

\begin{lem}\label{imbr}
Let $\A$ be an additive category with split idempotents and
products, which possesses a product generator $G$.
Denote $\CS=\{G[n]\mid n\in\Z\}$ the closure of $G$ under
suspensions and desuspensions in $\Htp\A$.
\begin{itemize}
\item[{\rm a)}] If given two composable maps $X\to Y\to Z$ whose
composition is $0$ in $\A$, then $X\to Y$ factors through a
subobject $Y'\leq Y$ such that the composed map $Y'\to Y\to Z$
vanishes, then $\Htp\A$ is $\CS$-cofiltered.

\item[{\rm b)}] If $\A$ has images or kernels, then $\Htp\A$
is $\CS$-cofiltered.
\end{itemize}
\end{lem}

\begin{pf} a) We will show inductively that a bounded complex
with less than $n+1$ non--zero entries is in $\Prd_n(\CS)$, where $n$ runs over all
positive integers. This is clear
for $n=0$, since $G$ is a product generator of $\A$. Now we suppose the
property true for any complex with
$\leq n$ non--zero entries. Let \[\cdots\to0\to X^0\to\cdots\to X^n\to0\to\cdots\]
be a bounded complex. The diagram \[\diagram
\cdots\rto&0\rto\dto&0\rto\dto&\cdots\rto&0\rto\dto&X^n\rto\dto^{=}&0\rto\dto&\cdots\\
\cdots\rto&0\rto\dto&X^0\rto\dto^{=}&\cdots\rto&X^{n-1}\rto\dto^{=}&X^n\dto\rto&0\rto\dto&\cdots\\
\cdots\rto&0\rto&X_0\rto&\cdots\rto&X^{n-1}\rto&0\rto&0\rto&\cdots
\enddiagram\]
is an exact sequence of complexes which splits in each degree.
According to \cite[Example 6.1]{Kel} it leads to a triangle
proving the induction step.

Finally consider an infinite complex \[X=\cdots\la
X^{n-1}\stackrel{d^{n-1}}\la X^n\stackrel{d^n}\la
X^{n+1}\la\cdots.\] By hypothesis, the map $d^{n-1}$ factors
through a subobject $Y^n\leq X^n$, such that $Y^n\la
X^n\stackrel{d^n}\la X^{n+1}$ vanishes, for all $n\in\Z$. For all
$i\in\N$, consider the bounded complex
\[X(i)=\cdots\to0\to Y^{-i}\to X^{-i}\to X^{-i+1}\to\cdots\to X^{i-1}\la
X^i\to0\to\cdots,\] and the map of complexes  $\epsilon(i):X(i+1)\to X(i)$
as in the following diagram:
\[\diagram \cdots\rto&0\rto&Y^{-i}\rto&X^{-i}\rto&\cdots\rto&X^i\rto&0\rto&\cdots\\
\cdots\rto&Y^{-i-1}\rto\uto&X^{-i-1}\rto\uto&X^{-i}\rto\udouble&\cdots\rto&X^i\rto\udouble&X^{i+1}\rto\uto&\cdots
\enddiagram\]
Applying \cite[Lemma 2.6]{IK} we infer
that $X$ is isomorphic in $\Htp\A$ to the homotopy limit of a the
chain of bounded complexes \[\cdots\la X(2)\stackrel{\epsilon(1)}\la
X(1)\stackrel{\epsilon(0)}\la X(0),\] thus $X$ is $\CS$-cofiltered.

b) We apply a) with $Y^n=\ima d^{n-1}$ or $Y^n=\ker d^n$, for all
$n\in\Z$.
\end{pf}

\begin{thm}\label{kmodbr} Let $\A$ be an additive category with products and
split idempotents,
possessing also images or kernels. Then $\Htp\A\opp$ satisfies Brown representability
if and only if $\A$ has a product generator.
In particular, if $R$ is a ring then $\Htp{\Mod{R}}\opp$ satisfies Brown
representability if and only if $\Mod{R}$ has a product generator.
\end{thm}

\begin{pf}
The direct implication is \cite[Theorem 2]{MS}, whereas the
converse follows by Lemma \ref{imbr} b) and Theorem \ref{cof}.
Finally note that the category
$\Mod{R}$ is additive with products and has both images and kernels.
\end{pf}

\begin{rem}
{\rm If the ring $R$ is pure--semisimple, then $\Mod R=\Add{G}$ for
some $G\in\Mod R$ (in fact $G$ is the direct sum of a family of
representatives of all isomorphism classes of finitely presentable
modules). In this case, $\Add{G}$ is closed under products, so $G$
is product--complete hence $\Add{G}=\Prod{G}$ (see \cite[Theorem
6.7]{KH}). Consequently $\Htp{\Mod{R}}\opp$ satisfies Brown
representability, by Theorem above. This was already known
since $\Mod{R}$ is a pure--semisimple finitely presentable
category which is closed under products, so it is compactly
generated by \cite[Theorem 5.2]{S}. It would be therefore
interesting to characterize the class of rings $R$ for which the module category
$\Mod R$ has a product generator. If we could indicate a non
pure--semisimple ring belonging to this class, then we would
produce an example of a triangulated category with products and
coproducts, namely $\K=\Htp{\Mod{R}}$ such that $\K\opp$, but not $\K$, satisfies
Brown representability. To the best of our knowledge,
such an example is yet unknown. Note added in proof: It seems that a ring $R$ for which 
$\Mod R$ has a product generator is pure--semisimple (see \cite{Br}), therefore Brown representability 
and its dual are equivalent for $\Htp{\Mod R}$.}
\end{rem}

\begin{rem}\label{kop}
 {\rm There is an isomorphism of categories $\Htp\A\opp\stackrel{\sim}\longrightarrow\Htp{\A\opp}$,
 which is easy to establish (for example, this is written down in \cite[Theorem 2.1.1]{MDer}). Applying
 this isomorphism of categories, we may dualize all results in this section. Thus we may conclude
 that if $\A$ is an additive category with split idempotents and coproducts,
 possessing also images or cokernels, then
 $\Htp\A$ satisfies Brown representability theorem if and only if $\A$ is pure--semisimple.
 Note that this statement
 is already known for $\A=\Mod{R}$, or more generally for a finitely accessible category
 with coproducts $\A$, as we may see by a combination between \cite[Theorem 1]{MS} and
 \cite[Proposition 2.6]{S}.
 However the results in \cite{MS} and \cite{S} may not be dualized in order to obtain
 Theorem \ref{kmodbr}, since the argument used there for showing that $\Htp\A$ satisfies
 Brown representability,
 where $\A$ is a
 pure--semisimple, finitely accessible additive category with coproducts goes as follows:
 If $\A$ enjoys
 all these properties, then $\Htp\A$ is well generated by \cite[Theorm 5.2]{S}, therefore
 it satisfies Brown
 representability by \cite[Theorem 8.3.3 and proposition 8.4.2]{N}.
 But none of the notions ``module category'',
 ``finitely accessible category'' and ``well generated triangulated category'' is self--dual.}
\end{rem}

\begin{rem}
 {\rm Let $R$ be a ring with $\gldim R\leq1$. Then the category $\Inj R$ of all injective modules
 is additive, closed under products, idempotents and images and every injective cogenerator of
 $\Mod R$ is a product generator for $\Inj R$. Thus Theorem \ref{kmodbr} gives another proof for
 the fact that $\Htp{\Inj R}\opp$ satisfies Brown representability. This fact is already known since
 $\Htp{\Inj R}$ is equivalent to the derived category which is compactly generated.}
\end{rem}

\begin{expl}
{\rm In the Introduction we said that this paper completes the picture in \cite{MP}.
Note that \cite[Theorem 3]{MP} gives an example of a triangulated coproduct
preserving functor which has no right adjoint, namely the inclusion functor
$\Htp\A\to\Htp\Ab$, where $\A$ is the full subcategory of all flat Mittag--Leffler
abelian groups.
Using the equivalence of categories
$\Htp\A\opp\stackrel{\sim}\longrightarrow\Htp{\A\opp}$ from Remark \ref{kop},
we obtain a triangulated product preserving functor which has no left adjoint.

Here we will provide another example of this kind, which holds only in an extension of
ZFC. More precisely, assume there are no measurable cardinals.
 For every cardinal $\lambda$ let us denote by $\Z^\lambda$ the product of $\lambda$-copies of $\Z$
and by $\Z^{<\lambda}$ its subgroup consisting of sequences with support (i.e. the set of non-zero entries) of
cardinality smaller then $\lambda$.
Let $\A\subseteq\Ab$ be the closure under products and direct factors of the class of all
abelian groups of the form $\Z^\lambda/\Z^{<\lambda}$, where $\lambda$ runs over all regular
cardinals.  The inclusion functor $\Htp\A\to\Htp\Ab$ is triangulated and preserves products.
If we suppose that it has a left adjoint then $\Htp\A$ must be preenveloping in
$\Htp\Ab$ by
Corollary \ref{envelops}. For $A\in\Ab$, the complex having $X$ in degree 0 and 0 elsewhere
must have an $\Htp\A$-preenvelope,
which is a complex $X$ with entries in $\A$. It is not hard to see that $X\to X^0$ is an
$\A$-preenvelope on $A$. But this
contradicts  \cite[Proposition 2.5]{CGR}, where it is shown that, under the hyopthesis of
nonexistence of measurable cardinals,
the class $\A$ is not preenveloping in $\Ab$.}

\end{expl}

% ------------------------------------------------------------------------------


\begin{thebibliography}{99}

\bibitem{B} A. Beligiannis, {Relative homological algebra and
    purity in triangulated categories\/}. {J. Algebra}, {\bf227} (2000), 268--361.
    
\bibitem{Br} S. Breaz, {A note on $\Sigma$-pure injective modules}, preprint arXiv:1304.0979v2 [math.RA].

\bibitem{CGR} C. Casacuberta, J. Guti\'erez, J. Rosick\'y, {Are all localizing subcategories of stable
homotopy category coreflective?}, preprint, arXiv:11062218v1 [math.CT].

\bibitem{CN} C. Casacuberta, A. Neeman, {Brown representability does not come for
free}, Math. Res. Lett., {\bf16}(2009), 1--6.

\bibitem{GT} R. G\"obel, J. Trlifaj, {\em Approximatins and Endomorphisms Algebras of Modules},
{De Gruyter Expositions in Mathematics}, {41}, Walter de Gruyter, 2006.


\bibitem{H} A. Heller, {On the representability of homotopy functors},
{J. London Math. Soc.}, {\bf23}, (1981), 551--562.

\bibitem{IK} S. Iyengar, H. Krause, {Acyclicity versus total acyclicity for complexes over noetherian rings},
{Doc. Math.} {\bf11} (2006), 207--240.

\bibitem{Kel} B. Keller, {Derived categories and their uses}, Chapther of {\em Handbook of Algebra},
edited by M. Hazewinkel, Elsevier 1996.

\bibitem{KH} H. Krause, {\em The spectrum of a module category}, Mem. Amer. Math. Soc.
{\bf149} (2001), x+125 pp.

\bibitem{KS} H. Krause, {Smashing subcategories and the telescope conjecture
-- an algebraic approach}, {Invent. Math.} {\bf139} (2000),
99--133.


\bibitem{KL} H. Krause, {Localization theory for triangulated
categories}, in  {\em Triangulated categories}, 161--235, London Math. Soc. Lecture Note Ser. 375, 
Cambridge Univ. Press, Cambridge, 2010.

\bibitem{MW} S. Maclane, {\em Categories for the Working
Mathematician}, Springer Verlag, New York, Heidelberg, Berlin,
1971.

\bibitem{MDer} D. Mili\v ci\'c, {Lectures on Derived Categories}, preprint available at author's homepage:
{\tt http://www.math.utah.edu/~milicic/}

\bibitem{MP} G. C. Modoi, {On perfectly generating projective classes in
triangulated categories}, Comm. Algeba, {\bf38}(2010), 995--1011.

\bibitem{MS} G. C. Modoi, J. \v S\v tov\'\i \v cek,
{Brown representability often fails for homotopy categories of
complexes}, {J. K-Theory}, {\bf9} (2012), 151--160.

\bibitem{N} A. Neeman, {\em Triangulated Categories\/}, Annals of
  Mathematics Studies, {\bf148}, Princeton University Press,
  Princeton, NJ, 2001.

\bibitem{NR} A. Neeman, {Brown Representability follows from
Rosick\'y's theorem}, {J. Topology}, {\bf2}(2009), 262--276.

\bibitem{NA} A. Neeman, {Some adjoints in homotopy categories}, {Ann. Math.}, {\bf171} (2010),
2143--2155.

\bibitem{P} M. Porta, {The Popescu--Gabriel theorem for triangulated categories},
{Adv. Math.}, {\bf225} (2010), 1669--1715.

\bibitem{S} J. \v S\v tov\'\i \v cek, {Locally well generated homotopy categories of
complexes}, Doc. Math. {\bf15} (2010), 507--525.
\end{thebibliography}
\end{document}